\documentclass[12pt]{article}
\usepackage{amsthm}
\usepackage{graphicx}

\newtheorem{theorem}{Theorem}
\newtheorem{lemma}[theorem]{Lemma}
\newtheorem{corollary}[theorem]{Corollary}

\title{Continued fractions and Catalan problems}
\author{Mahendra Jani \\William Paterson University \and Robert G. Rieper \\William Paterson University}

\begin{document}
\maketitle
% ----------------------------------------------------------------
\begin{abstract}
We find a generating function expressed as a continued fraction that enumerates ordered trees by the number of vertices
at different levels. Several Catalan problems are mapped to an ordered-tree problem and their generating functions also
expressed as a continued fraction. Among these problems is the enumeration of $(132)$-pattern avoiding permutations
that have a given number of increasing patterns of length $k$. This extends and illuminates a result of Robertson, Wilf
and Zeilberger for the case $k=3$.
\end{abstract}

\section{Introduction}

A Catalan problem is any enumerative problem that produces the Catalan sequence of numbers or one of its many
$q$-analogs. Stanley~\cite{RS99} provides a catalog of $66$ Catalan problems.  Interestingly, many of the generating
functions that arise from these problems can be given as a continued fraction with a simple yet elegant form.  Two of
these generating functions are reproduced below and a third we derive anew.  Our intent is to show that the first two
continued fractions are special instances of the third with the implication that many others are as well.  We begin
with the three Catalan problems and their corresponding generating functions.

\textbf{Problem 1}. A $(132)$ pattern (respectively, a $(123)$ pattern) in a permutation $\pi$ of length $n$ is a
triple $1\le i<j<k\le n$ of indices for which $\pi(i)<\pi(k)<\pi(j)$ (respectively, $\pi(i)<\pi(j)<\pi(k)$).  Let
$f_r(n)$ denote the number of permutations $\pi$ of length $n$ that have no $(132)$ patterns and exactly $r$ $(123)$
patterns. Recently, Robertson, Wilf and Zeilberger~(\cite{RWZ99}) derived the generating function

\begin{equation}
\label{e:gen1} \sum_{r,n\ge 0}f_r(n)z^nq^r ={1\over\displaystyle 1-{\strut z\over\displaystyle 1-{\strut
z\over\displaystyle 1-{\strut zq\over\displaystyle 1-{\strut zq^3\over\displaystyle 1-{\strut zq^6\over\displaystyle
\ddots}}}}}}
\end{equation} in which the $l$th numerator is $zq^{{l-1\choose 2}}$ ($l$ is for \emph{level} in anticipation of Problem 3).

\textbf{Problem 2}.  The number of lattice paths from $(0,0)$ to $(n,n)$ with steps $(1,0)$ and $(0,1)$ that never rise
above the line $y=x$ is a Catalan number.  Let $P$ be such a path, $A(P)$ the area under the path (and above the
$x$-axis), and let $C_n(q)=\sum_P q^{A(P)}$.  Then a generating function is given by (see Exercise 6.34 in
Stanley~\cite{RS99} and replace the $x$ therein with $zq$ )

\begin{equation}
\label{e:gen2}
 \sum_{n\ge 0}q^{n+1 \choose 2}C_n(1/q)z^n ={1\over\displaystyle 1-{\strut zq\over\displaystyle 1-{\strut
zq^2\over\displaystyle 1-{\strut zq^3\over\displaystyle 1-{\strut zq^4\over\displaystyle \ddots}}}}}
\end{equation} in which the $l$th numerator is $zq^l$.

\textbf{Problem 3}.  The number of ordered trees (also known as plane trees) on $n$ edges is a Catalan number.  The
\emph{level} of a vertex is the number of edges on the unique path from the root to the vertex.  Thus, the root is the
unique vertex at level zero and the vertices at level one are adjacent to the root.  Let $T_{l_1,l_2,\ldots}$ be the
number of ordered trees that have $l_k$ vertices at level $k> 0$ and let $v_k$, $k>0$, be indeterminates.  The
generating function $T$ that enumerates ordered trees by the number of vertices at each level is defined as

$$ T(v_1,v_2,\ldots )=\sum_{l_1,l_2,\ldots\ge 0}T_{l_1,l_2,\ldots}v_1^{l_1}v_2^{l_2}\cdots . $$ The first few terms
(number of edges $n\le 3$) of $T$ are $$ T(v_1,v_2,\ldots )=1
+v_1+v_1v_2+v_1^2+v_1v_2v_3+2v_1^2v_2+v_1v_2^2+v_1^3+\cdots . $$ That $T$ can be written as a continued fraction that
subsumes the previous continued fractions is our main result.  It is simple, yet has some interesting applications.

\begin{theorem}
\label{t:main} The generating function that enumerates ordered trees by the number of vertices at each level is $$
T(v_1,v_2,\ldots )= {1\over\displaystyle 1-{\strut v_1\over\displaystyle 1-{\strut v_2\over\displaystyle 1-{\strut
v_3\over\displaystyle  \ddots}}}}$$
\end{theorem}

\begin{proof}
We exploit the natural recursive property of ordered trees to obtain a recursion for $T$.  The recursion immediately
leads to the continued fraction.  Any ordered tree on more than one vertex can be constructed from a collection of
others (the subtrees) by joining the roots of these subtrees to a new vertex.  The new vertex becomes the root of the
tree under construction.  Note that the level of a vertex in a subtree \emph{increases by one} after the new root is
inserted. The function $T(v_1,v_2,\ldots )$ enumerates the choices for a subtree and each of these choices contributes
a factor of $v_1T(v_2,v_3,\ldots )$ because of the level changes.  The factor of $v_1$ is present because the root of
the subtree becomes a vertex at level one.  Thus, the trees with $k$ subtrees (of the root) are enumerated by
$v_1^kT^k(v_2,v_3,\ldots )$. The generating function satisfies
\begin{eqnarray*}
T(v_1,v_2,\ldots )  &=& 1+v_1T(v_2,v_3,\ldots )+v_1^2T^2(v_2,v_3,\ldots )+\ldots \\
                        &=&\frac{1}{1-v_1T(v_2,v_3,\ldots )}.
\end{eqnarray*}
Iteration of the last functional recursion produces the continued fraction.
\end{proof}

An immediate application is obtained by replacing each indeterminate $v_k$ with $z$ denoted simply as $T(z)$. The
resulting function enumerates ordered trees by the number of edges and is
\begin{eqnarray*}
T(z )   &=& {1\over\displaystyle 1-{\strut z\over\displaystyle 1-{\strut z\over\displaystyle 1-{\strut
z\over\displaystyle \ddots}}}} \\
        &=& \frac{1}{1-zT(z)}.
\end{eqnarray*}
The well-known solution of the above generates the Catalan numbers and is $T(z)=(1-\sqrt{1-4z})/2z$.

The more challenging applications are the evaluations needed to produce the continued fractions of permutations
(Equation~\ref{e:gen1} of Problem 1) and lattice paths (Equation~\ref{e:gen2} of Problem 2).  Both applications require
that we map their respective problems to an ordered-tree problem. These mappings are of interest in their own right and
we explore them now.  We begin with the lattice path problem because it is simpler and the mapping is already known.

% section2.tex

\section{Lattice paths and ordered trees}
We draw our ordered trees with the root at the top and proceed downward to the leaves.  The first leaf is the leftmost
leaf in the drawing and the remaining leaves are referred to by their positions in a left-to-right order.  A preorder
(depth-first) traversal of the ordered tree provides a well-known correspondence with a lattice path. When an edge of
the tree is traversed downward away from the root we take a $(1,0)$ (east) step in the lattice, otherwise we take a
$(0,1)$ (north) step.  In this manner, an ordered tree with $n$ edges corresponds to a unique lattice path from $(0,0)$
to $(n,n)$.

If the path $P$ corresponds to the tree $T$, then we need to determine what statistic of the tree corresponds to the
area $A(P)$ under the path.  We let $A(T)=A(P)$ be this statistic of the tree and claim that it depends only on the
 vertex levels.

\begin{lemma}
\label{l:paths} If $T$ is an ordered tree on $n$ edges, then $$A(T)={n+1 \choose 2}-\sum_{\mathrm{vertices }\,
v}\mathrm{level}(v),$$ where $\mathrm{level}(v)$ is the level of vertex $v$.
\end{lemma}

\begin{proof}
Let $w$ be the rightmost leaf of the tree $T$.  Our immediate interest is to calculate the area under the east step
that arises in the lattice path by traversing the last edge to this leaf downward away from the root.  This area is
equal to the height that the east step has in the lattice path and is equal to the number of north steps that have
occurred prior to the east step. There are a total of $n$ north steps in the lattice path (one for each edge of the
tree) and $\mathrm{level}(w)$ remaining north steps after the east step.  Hence, the east step is at height
$n-\mathrm{level}(w)$ in the lattice. If we now delete the leaf $w$ from the tree $T$, then the resulting tree has a
lattice path with area $n-\mathrm{level}(w)$ less than that of $T$.  A formal inductive argument on the number of edges
provides the result.
\end{proof}

We use the lemma to prove the following continued-fraction result.  The result is the same as that given in
Equation~\ref{e:gen2}.
\begin{theorem}
\label{t:area} If $C_n(q)=\sum_Tq^{A(T)}$ enumerates the set of ordered trees on $n$ edges by the area under their
corresponding lattice paths, then $$\sum_{n\ge 0}q^{n+1 \choose 2}C_n(1/q)z^n ={1\over\displaystyle 1-{\strut
zq\over\displaystyle 1-{\strut zq^2\over\displaystyle 1-{\strut zq^3\over\displaystyle 1-{\strut zq^4\over\displaystyle
\ddots}}}}}$$ in which the $l$th numerator is $zq^l$.
\end{theorem}

\begin{proof}
Let $T$ be an ordered tree on $n$ edges and assign to every nonroot vertex at level $l>0$ the value $zq^l$.  The
product of all these values is then $z^nq^{\sum_v\mathrm{level}(v)}$. Summing over all ordered trees on $n$ edges we
have by the lemma $$ z^n\sum_T q^{\sum_v\mathrm{level}(v)}=z^n\sum_Tq^{{n+1 \choose 2}-A(T)}=z^nq^{n+1 \choose
2}C_n(1/q).$$ The sum of these over all $n\ge 0$ then enumerates ordered trees and the generating function is given by
the continued fraction of Theorem~\ref{t:main} with $v_l=zq^l$.
\end{proof}

%  section3.tex

\section{Permutations and ordered trees}
  The previous problem used an existing bijection between the set of ordered trees and the set of lattice paths to get
  the desired result.  We seek a similar approach for ordered trees and permutations.  There are many ways to map a
  permutation onto a tree (often an unordered tree) but none of these serve our needs.  The mapping we introduce
  appears to be new.

\begin{figure}
\begin{center}
\includegraphics[]{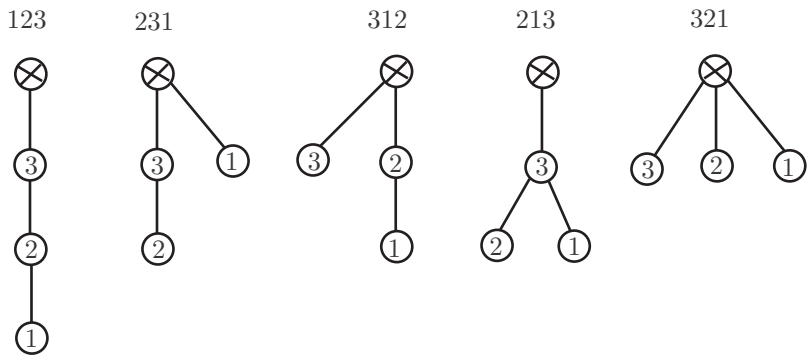}
\caption{The five ordered trees on $3$ edges and their permutations.}\label{f:one}
\end{center}
\end{figure}

  Let $T$ be an ordered tree on $n$ edges.  We use a preorder traversal of $T$ to label the nonroot vertices in
  decreasing order with the integers $n,n-1,\ldots,1$.  Thus, the first vertex visited gets the label $n$ and the last
  receives $1$.  We now construct a permutation written as a word by \emph{reading} the labeled tree in
  \emph{postorder}.  We again traverse the tree from left to right and record the label of a vertex when we last visit
  it.  In Catalan fashion, the five ordered trees and their corresponding permutations are shown in Figure~\ref{f:one}.
  Note that the only permutation missing from those of length three is $132$.  The (132) pattern has been avoided. Also
  note that there is exactly one permutation with a (123) pattern (the first permutation shown).  Thus, recalling the
  definition of $f_r(n)$ given in the introduction we have $f_0(3)=4$, $f_1(3)=1$, and $f_r(3)=0,r>1$. We generalize
  these observations after introducing some useful notation.

If $T$ is an ordered tree on $n$ edges, then we let $\pi(T)$ be its corresponding permutation written as a word on the
numbers $1,2,\ldots,n$.  We let $\pi(T,k)$ be the same permutation except we use the corresponding numbers
$1+k,2+k,\ldots,n+k$. For example, the first tree shown in the figure has $\pi(T)=\pi(T,0)=123$ and $\pi(T,1)=234$. For
emphasis, we denote the concatenation of two words $\pi$ and $\pi'$ as $\pi\wedge\pi'$ (instead of the usual
$\pi\pi'$).  The following lemma describes how the permutation of a tree can be constructed from those of its subtrees.

\begin{lemma}
\label{l:perm} Let $T$ be an ordered tree on $n$ edges with subtrees $T_1,T_2,\ldots, T_s$ on $n_1,n_2,\ldots,n_s$
edges, respectively.  Let $N_0=n$ and $N_k=N_{k-1}-n_k-1,k=1,2,\ldots,s$, then $$\pi(T)=\pi(T_1,N_1)\wedge
N_0\wedge\pi(T_2,N_2)\wedge N_1\wedge\ldots\wedge\pi(T_s,N_s)\wedge N_{s-1}. $$
\end{lemma}

\begin{proof}
Note that $N_0=n$ is the total number of vertices to be labeled, that $N_1$ is the total number of vertices to be
labeled after those of the subtree $T_1$, and so on.  Since the vertices of $T$ are labeled in decreasing order using a
preorder traversal, $i<j$ implies that all vertices of $T_i$ receive labels greater than those of $T_j$. Since $\pi(T)$
is constructed by reading these labels in postorder, $i<j$ implies that the vertex labels in $T_i$ appear in $\pi(T)$
before any of those in $T_j$. Thus, $\pi(T)$ begins as $\pi(T_1)$ with each number incremented by $N_1$, i.e., as
$\pi(T_1,N_1)$. The root of this subtree is the first vertex visited in the preorder traversal and receives the label
$n=N_0$.  It is read last among the vertices of $T_1$ in postorder, however, so that $\pi(T)$ begins as
$\pi(T_1,N_1)\wedge N_0$. The general case follows similarly.
\end{proof}

Note that the lemma also provides the means to prove that the mapping $T\rightarrow \pi(T)$ is injective.  That it is
onto the set of all $(132)$-pattern avoiding permutations is proved in the next theorem.  Before proceeding to this
theorem we present another lemma which enables us to count $(123)$ patterns and their generalization. An
\emph{increasing} pattern of \emph{length} $k$, $k>0$, in a permutation $\pi$ of length $n$ is a $k$-tuple  $1\le
i_1<i_2<\ldots <i_k\le n$ of indices for which $\pi(i_1)<\pi(i_2)<\ldots <\pi(i_k)$.

\begin{lemma}
\label{l:increase} Let $T$ be an ordered tree on $n$ edges and $V_k$ a subset of $k$ vertices, $0<k\le n$, none of
which are the root, then the labels of these vertices provide an increasing pattern of length $k$ in $\pi(T)$ if and
only if they lie along a path from the root to some leaf.
\end{lemma}

\begin{proof}
We induct on $n$.  If $n=1$, then $k=1$ and the lemma is obvious.  In fact, this is the case for all $n$ whenever
$k=1$. We assume that the lemma is true for any ordered tree on $n$ or fewer edges, $n>0$, and let $T$ be a tree on
$n+1$ edges. Suppose that $\pi(T)$ contains an increasing pattern of length $k>1$. We let $v_i$ be the vertex in $T$
that provides the $i$ in such a pattern, $i=1,2,\ldots,k$. Then $v_k$ must receive a larger label than $v_1$ implies
that the subtree containing $v_k$ precedes or is the same as that of $v_1$ (see the proof of Lemma~\ref{l:perm}).
However, the label of $v_k$ must be read after that of $v_1$ implies that its subtree must follow or be the same as
that of $v_1$. The conclusion is that they are in the same subtree.  A similar argument applies to $v_i$ and $v_1$,
$i>1$, so that all the vertices must be in the same subtree, say $T_i$. Thus, the increasing pattern of length $k$ lies
entirely within $\pi(T_i,N_i)\wedge N_{i-1}$.

Recall that it is the root of the subtree that receives the label $N_{i-1}$.  If the root of the subtree is one of the
vertices providing the increasing pattern, then it must be $v_k$. We must consider two cases depending on whether the
subtree root is $v_k$ or not.

If it is not, then the pattern lies entirely within $\pi(T_i,N_i)$ and corresponds uniquely to an increasing pattern of
length $k$ in $\pi(T_i)$.  By the inductive hypothesis, the vertices $v_1,v_2,\ldots,v_k$ must lie along some path from
the root of the subtree to a leaf. Necessarily, this is also a path from the root of $T$ to a leaf as required by this
lemma.

If the root of the subtree is $v_k$, then $v_1,v_2,\ldots,v_{k-1}$ provide the labels for an increasing pattern of
length $k-1$ in $\pi(T_i,N_i)$. This pattern corresponds to a unique increasing pattern of length $k-1$ in $\pi(T_i)$
and again by hypothesis the vertices $v_1,v_2,\ldots,v_{k-1}$ must lie along a path from the root of the subtree to
some leaf of the subtree.  This path together with the root of the subtree provides the needed path of this lemma.

Conversely, if $V_k=\{v_1,v_2,\ldots,v_k\}$ is a subset of nonroot vertices of $T$, where we may assume the label of
$v_i$ is less than that of $v_j$ whenever $i<j$, and the vertices lie along a path from the root to a leaf, then $v_i$
lies below $v_j$ on this path whenever $i<j$.  Thus, when the labels are read in postorder, the label of $v_i$ is read
prior to that of $v_j$, $i<j$.  The vertices then provide an increasing pattern of length $k$ in $\pi(T)$.

\end{proof}

The two previous lemmas enable us to prove the following interesting combinatorial theorem.  Its corollary establishes
a continued fraction as the generating function of $(132)$-avoiding permutations by number of increasing patterns of
length $k$.

\begin{theorem}
\label{t:perm} A permutation $\pi$ avoids the $(132)$ pattern if and only if $\pi=\pi(T)$ for some tree $T$.  If this
is the case, then the number of increasing patterns of length $k$ depends only on the levels of the vertices in the
tree and is given by $\sum_v {\mathrm{level}(v)-1 \choose k-1}$.
\end{theorem}

\begin{proof}
Suppose that $\pi(T)$ contains a $(132)$ pattern and that $T$ is among the smallest such trees. We let $v_3$ be a
vertex in $T$ that provides the $3$ in such a pattern and let $v_1,v_2$ be the vertices that provide the corresponding
$1$ and $2$, respectively. Then $v_3$ must receive a larger label than $v_1$ implies that the subtree containing $v_3$
precedes or is the same as that of $v_1$ (see the proof of Lemma~\ref{l:perm}).  However, the label of $v_3$ must be
read after that of $v_1$ implies that its subtree must follow or be the same as that of $v_1$. The conclusion is that
they are in the same subtree.  A similar argument applies to $v_2$ and $v_1$ so that all three vertices must be in the
same subtree, say $T_i$.

Also note that none of them can be the root of the subtree since the root receives the largest label among the vertices
of the tree and, hence, can not be $v_1$ or $v_2$. Its label appears later in $\pi(T)$ than the others in the subtree
implies that the root can not be $v_3$.  Thus, the $(132)$ pattern lies entirely within $\pi(T_i,N_i)$ which implies
that $\pi(T_i)$ itself must have a $(132)$ pattern contradicting our choice of $T$. Since the number of $(132)$-pattern
avoiding permutations of length $n$ and the number of ordered trees on $n$ edges are the same Catalan number, the
mapping $T\rightarrow\pi(T)$ is a bijection between these sets. It is an instructive exercise to construct $T$ from a
$(132)$-pattern avoiding permutation.

It remains to determine the number of increasing patterns of length $k$ in $\pi (T)$.  As a result of
Lemma~\ref{l:increase} it is only necessary to count the number of vertex subsets of size $k$, none of which are the
root, such that the vertices lie along a path from the root to a leaf.  We claim this number is $\sum_v
{\mathrm{level}(v)-1 \choose k-1}$ as stated. To see this, let $v$ be any vertex of $T$ and choose $v_k=v$.  There are
$\mathrm{level}(v)-1$ nonroot vertices other than $v$ along the unique path from the root to $v$.  From these we select
any $k-1$ of them, which together with $v$, form the required subset.  It is clear that every subset with the required
properties arises this way and we are done.

\end{proof}

We now use the theorem to write a generating function as a continued fraction.  We let $f_r^{(k)}(n)$ denote the number
of permutations of length $n$ that have no $(132)$ pattern and exactly $r$ increasing patterns of length $k$. The case
$k=3$ is that considered by Robertson, Wilf and Zeilberger~\cite{RWZ99}.

\begin{corollary}
\label{c:perm} The generating function that enumerates $(132)$-pattern avoiding permutations of length $n$ by number of
increasing patterns of length $k$ is
\begin{equation}
\label{e:gen3} \sum_{r,n\ge 0}f_r^{(k)}(n)z^nq^r ={1\over\displaystyle 1-{\strut N_1\over\displaystyle 1-{\strut
N_2\over\displaystyle 1-{\strut N_3\over\displaystyle 1-{\strut N_4\over\displaystyle \ddots}}}}}
\end{equation} in which the $l$th numerator $N_l$ is $zq^{{l-1\choose k-1}}$.
\end{corollary}

\begin{proof}
Let $T$ be an ordered tree on $n$ edges and assign to every nonroot vertex at level $l>0$ the value $zq^{{l-1 \choose
k-1} }$.  The product of all these values is then $z^nq^{\sum_v{\mathrm{level}(v)-1 \choose k-1}}$ and the result
follows from the previous theorem and Theorem~\ref{t:main}.
\end{proof}

% ----------------------------------------------------------------
\providecommand{\bysame}{\leavevmode\hbox to3em{\hrulefill}\thinspace}

% ----------------------------------------------------------------

\vspace{.5in} \noindent \texttt{janim@wpunj.edu}\\ \texttt{jrieper@cybernex.net}\\
\end{document}